\documentclass[11pt,a4paper,final]{amsart}

\usepackage{url} 
\usepackage[hidelinks]{hyperref}
\usepackage[latin1]{inputenc}
\usepackage{amssymb,amsthm,amsmath,amstext,amsfonts,mathtools}
\usepackage{amsaddr}
\usepackage{graphicx}
\usepackage[a4paper,twoside=true,top=2.0cm,bottom=2.2cm,left=2.5cm,right=2.5cm]{geometry}
\usepackage[normalem]{ulem}

\def\Rset{{\mathbb{R}}}
\def\Cset{{\mathbb{C}}}
\def\Zset{{\mathbb{Z}}}
\def\Nset{{\mathbb{N}}}
\def\eu{\ensuremath{\mathrm{e}}}
\def\iu{\ensuremath{\mathrm{i}}}
\def\du{\ensuremath{\mathrm{d}}}

\newtheorem{theorem}{Theorem}
\newtheorem{remark}[theorem]{Remark}
\newtheorem{proposition}[theorem]{Proposition}
\newtheorem{lemma}[theorem]{Lemma}
\newtheorem{corollary}[theorem]{Corollary}
\theoremstyle{definition}

\def\texttiny#1{{\text{\tiny{#1}}}}
\def\DC{{}^{\texttiny{C}}\! D}

\begin{document}

\title[On initial conditions for fractional delay differential equations] 
{On initial conditions for \\ fractional delay differential equations}

%
\author{Roberto Garrappa}
\address[Roberto Garrappa]{Dipartimento di Matematica, Universit\`a degli Studi di Bari\\ Via E. Orabona 4 - 70126 Bari - Italy\\Member of the INdAM Research group GNCS}
\email[R.~Garrappa]{roberto.garrappa@uniba.it}

\author{Eva Kaslik}
\address[Eva Kaslik]{Department of Mathematics and Computer Science West University of Timi\c{s}oara \\ Bd. V. P\^{a}rvan 4, 300223 Timi\c{s}oara, Romania}
\email[E.~Kaslik]{eva.kaslik@e-uvt.ro}

\thanks{This is the preprint of a paper published in [Commun. Nonlinear Sci. Numer. Simul., (2020)] and available at \url{https://doi.org/10.1016/j.cnsns.2020.105359}. This research was funded by the COST Action CA 15225 - ``Fractional-order systems- analysis, synthesis and their importance for future design''. The work of R.Garrappa is also partially supported by a GNCS-INdAM 2020 Project.}




\maketitle

\begin{abstract}
Derivatives of fractional order are introduced in different ways: as left-inverse of the fractional integral or by generalizing the limit of the difference quotient defining integer-order derivatives. Although the two approaches lead (under standard smoothness assumptions) to equivalent operators, the first one does not involve the function at the left of the initial point where, instead, the latter forces the function to assume selected values. With fractional delay differential equations new problems arise: the presence of the delay imposes to assign the solution not just at the initial point but on an entire interval. Due to the freedom in the choice of the initial function, some inconsistencies with the values forced by the fractional derivative are possible and the operators may no longer be equivalent. In this paper we discuss the initialization of fractional delay differential equations, we investigate the effects of the initial condition not only on the solution but also on the fractional operator as well and we study the difference between solutions obtained by incorporating or not the initial function in the memory of the fractional derivative. The exact solution of a family of linear equations is obtained by the Laplace transform whilst numerical methods are used to solve nonlinear problems; the different results are therefore shown and commented.
\end{abstract}


\def\Rset{{\mathbb{R}}}
\def\Cset{{\mathbb{C}}}
\def\Zset{{\mathbb{Z}}}
\def\Nset{{\mathbb{N}}}
\def\eu{{\ensuremath{\mathrm{e}}}}
\def\iu{{\ensuremath{\mathrm{i}}}}
\def\du{{\ensuremath{\mathrm{d}}}}
\def\texttiny#1{{\text{\tiny{#1}}}}
\def\DC{{}^{\texttiny{C}}\! D}
\def\DR{{}^{\texttiny{RL}}\! D}
\def\DG{{}^{\texttiny{GL}}\! D}
\def\DL{{}^{\texttiny{L}}\! D}
\def\phitau{{\phi_{\tau}}}
\def\Dphi{{}^{{\phi_{\tau}}}\! D}
\def\fR{{f}^{\text{{\scriptsize{R}}}}} 
\def\fC{{f}^{\text{{\scriptsize{C}}}}} 
\def\yR{{y}^{\text{{\scriptsize{R}}}}} 
\def\yC{{y}^{\text{{\scriptsize{C}}}}} 
\def\Cc{{}^{{\phi}}\! C_{-\tau}}
\def\yphi{{y}^{{{\phi_{\tau}}}}}

\section{Introduction}

In the last decades the interest toward models incorporating derivatives of fractional (i.e., non-integer) order has increased in a noteworthy way. Indeed, fractional differential equations (FDEs) describe in a more accurate way anomalous relaxation processes in which the external excitation does not have an instantaneous effect but depends on the past history of the system. Fields in which FDEs are satisfactorily employed range from biology to control theory, engineering, finance, optics, physics and so on.  

In more complex interconnected systems, time delays are also introduced since changes in one variable may affect other variables with some lags. For example, in feedback control systems, time delays account for delayed feedback. In models arising from biology, delays are incorporated to describe incubation time or time to maturity. In economic models, delays are included for a better modeling of transportation or information transmission.


Unlike standard integer-order delay differential equations (DDEs) for which well-grounded theories have been already established (see, for instance, \cite{BellenZennaro2003,BredaMasetVermiglio2015,Driver1977,Kuang1993} and references therein), the analysis and application of models incorporating fractional delay differential equations (FDDEs) rely on a theoretical basis which still presents some not completely clear aspects.

This is the case of the initial conditions to couple to the selected fractional derivative. Different definitions are indeed possible and each of them demands for specific conditions to initialize the problem. For instance, while the Caputo derivative allows classical initial conditions of Cauchy type, FDEs with the Riemann-Liouville (RL) derivative require initial conditions expressed as fractional derivatives of the unknown function. A further problem however arises with FDDEs. The presence of a delay demands not just a finite number of conditions at the initial point but it is also necessary to assign an initial function over an interval to the left of the initial point. 


A generic FDDE can be formulated as
\begin{equation}\label{eq:FDDE_General}
\begin{cases}
	\DC^{\alpha}_{0} y(t) = g(t,y(t),y(t-\tau)) &,~ t>0 \\
	y(t) = \phi(t) & ,~ -\tau \le t \le 0
\end{cases}
\end{equation}
where $\tau>0$ is a constant delay and $g:[0,T]\times \Rset \times \Rset \to \Rset$ and $\phi(t):[-\tau,0]\to\Rset$ are given functions. Just for simplicity we focus here on FDDEs (\ref{eq:FDDE_General}) with the Caputo fractional derivative $\DC^{\alpha}_{0}$ of order $0<\alpha<1$ but it is possible to extend the discussion to more involved problems (e.g., other fractional derivatives, higher order derivatives, more or non constant delays, so on).

There are different ways to derive the Caputo derivative. One approach introduces $\DC^{\alpha}_{0}$ as the left-inverse operator of the fractional RL integral and substantially ignores the behavior of the function to the left of the origin. An alternative approach, instead generalizes the limit of the difference quotient defining integer-order derivatives; this generalization, which is known as the Gr\"unwald-Letnikov (GL) fractional derivative, forces the function to assume suitable values (the constant $y_0=y(0)$ in our case) to the left of the origin to ensure the convergence of the series in the limit defining the GL derivative.


Under standard continuity assumptions, the two approaches lead to equivalent operators. However, whenever the function is forced to assume specific values in $[-\tau,0]$ due to the initial condition $\phi(t)$, the two approaches may no longer lead to equivalent operators. 

It is therefore necessary to analyze the impact on the derivative of the initial condition $y(t) = \phi(t)$, $t \in [-\tau , 0]$, and try to establish which operator is more appropriate in order to avoid inconsistencies. We think that clarifying this aspect may help other researchers in handling fractional derivatives in systems with delay.

Obviously, one could initialize the derivative at $-\tau$ instead of at 0, since the process described by (\ref{eq:FDDE_General}) actually originates at $-\tau$. However, changing the starting point modifies the derivative and leads to a different problem; we prefer here to just discuss the consequences of coupling $\DC^{\alpha}_{0}$ with the condition $y(t) = \phi(t)$, $t \in [-\tau, 0]$, as usually proposed in the literature, and refer to \cite{AcharLorenzoHartley2007} for results related to changes in the initial point.


We must also mention that the initialization of FDEs is, in general, an open and debated issue, earning considerable attention \cite{ortigueira2010system,mozyrska2011modified,ortigueira2012usefulness,sabatier2017misconceptions,maamri2017improved}. For instance, even though initial conditions associated to the RL derivative have previously been thought to lack physical meaning, Heymans and Podlubny \cite{heymans2006physical} contradicted this claim, by offering pertinent examples and by introducing the concept of ``inseparable twins''. On the other hand, Caputo's derivative has been criticized \cite{HartleyTrigeassouLorenzoMaamri2013} since its initial condition assumes non-vanishing values for all negative time, thus implicitly leading to systems which require an infinite amount of energy to initialize; with this respect, Lorenzo and Hartley have studied possible corrections to the initialization process (e.g., see \cite{LorenzoHartley2008}). Assigning proper initial conditions to fractional-order differential equations is of utmost importance with respect to the solution of the problem, as well as its interpretation in the framework of the modeled real world phenomenon. However, it is beyond the scope of this paper to discuss the general problem of the initialization of FDEs; we prefer to just focus on specific issues related to FDDEs.

The paper is organized in the following way. In Section \ref{S:DelayedFDEs} we recall some preliminary results from the theory of FDDEs. A general discussion regarding the way by which some of the most used fractional derivatives are obtained is presented in Section \ref{S:GeneralizingDerivatives}. We therefore observe that fractional derivatives force the function to assume values to the left of the initial point which may conflict with the initial function $\phi(t)$ which is, instead, a data from the model. To analyze the consequences of initial data conflicting with the values imposed by the fractional derivative, by means of the Laplace transform in Section \ref{S:ExactSolution} we first derive the exact solution of some linear test FDDEs with different initial conditions. The solutions of the same FDDEs are also obtained in Section \ref{S:GL_Modified} after modifying the fractional derivative in order to take into account the values of the initial function: the obtained results are therefore compared and a similar comparison is made for nonlinear problems by means of numerical methods in Section \ref{S:Nonlinear}. Some concluding remarks are finally presented in Section \ref{S:ConcludingRemarks} and an alternative proof concerning the exact solution of linear FDDEs, and based on induction arguments, is presented for completeness in an appendix at the end of the paper.

\section{Fractional delay differential equations}\label{S:DelayedFDEs}

Let $C([a,b],\Rset)$ denote the space of continuous functions $\xi:[a,b]\rightarrow \Rset$ endowed with the sup norm
\[\|\xi\|_\infty:=\sup_{t\in[a,b]}|\xi(t)|.\]

Throughout this paper we assume that in the initial value problem \eqref{eq:FDDE_General}, the function $g:[0,T]\times \Rset\times\Rset\rightarrow\Rset$ is continuous and $\phi\in C([-\tau,0],\Rset)$. 

A function $\varphi(\cdot,\phi)\in C([-\tau,T),\Rset)$ is a solution of the initial value problem \eqref{eq:FDDE_General} if 
\[
\begin{cases}
	\DC^{\alpha}_{0} \varphi(t,\phi) = g(t,\varphi(t,\phi),\varphi(t-\tau,\phi)) &,~ 0<t\leq T \\
	\varphi(t,\phi)= \phi(t) &,~-\tau \le t \le 0 \\
\end{cases}
\]

Several recent papers \cite{abbas2011existence,yang2013initial,morgado2013analysis,wang2016existence} have investigated the existence and uniqueness problem related to the initial value problem \eqref{eq:FDDE_General}. The most general result up to date regarding the existence and uniqueness of a global solution of the initial value problem \eqref{eq:FDDE_General} has been proved in \cite{cong2017existence}, subject to a  mild Lipschitz condition on the function $g$ with respect to its second (non-delay) variable.  Moreover, if the initial value problem \eqref{eq:FDDE_General} is considered on the whole semi-axis $[-\tau,\infty)$, the exponential boundedness of the global solution, which is mandatory when Laplace transform techniques are used in the qualitative analysis (as in the remainder of this paper), has also been proved. 

\begin{theorem}[see Theorem 4.1 in \cite{cong2017existence}]\label{thm.exp.bounded} If $g:[0,\infty)\times\Rset\times\Rset\rightarrow\Rset$ is continuous and satisfies the following two assumptions: 
\begin{itemize}
    \item[(H1)] There exists a positive constant $L>0$ such that 
    \[|g(t,x,y)-g(t,\hat{x},\hat{y})|\leq L(|x-\hat{x}|+|y-\hat{y}|),\quad\forall t\geq 0,~x,y,\hat{x},\hat{y}\in\Rset.\]
    \item[(H2)] There exists a constant $\beta>2L$  such that 
    \[\sup_{t\geq 0}\frac{\int_0^t (t-r)^{\alpha-1}|g(r,0,0)|dr}{E_\alpha(\beta t^\alpha)}<\infty.\]
\end{itemize}
Then the unique global solution $\varphi(\cdot,\phi)$ defined on the interval $[-\tau,\infty)$ of the initial value problem \eqref{eq:FDDE_General} is exponentially bounded; more precisely, there exists a constant $C>0$ such that
\[|\varphi(t,\phi)|\leq CE_\alpha(\beta t^\alpha),\quad\forall~ t\geq 0.\]
\end{theorem}

Nevertheless, another important problem which arises in the setting of FDDEs is the asymptotic behaviour of solutions. Stability and asymptotic properties of linear systems of fractional-order differential equations involving delayed and/or non-delayed terms have been explored in \cite{krol2011asymptotic,CermakHornicekKisela2016,cermak2019oscillatory}. Moreover, a linearization theorem has been proved in \cite{tuan2018linearized}, showing that an equilibrium of a nonlinear Caputo FDDE is asymptotically stable if
its linearization at the equilibrium is asymptotically stable.

Recent applications of FDDEs in the modelling of real world phenomena include the investigation of neural network models \cite{wang2017stability,huang2018effects}, gene regulatory networks \cite{tao2018hopf}, multi-agent systems  \cite{liu2019stability} epidemiological models \cite{latha2017fractional} and competition models \cite{leung2014periodic,wang2019stability}.

\section{Generalizing integer-order derivatives to fractional order}\label{S:GeneralizingDerivatives}

To focus on the main issues related to the initialization of FDDEs, we first present, in some detail, the mathematical process at the basis of the introduction of fractional derivatives. 

A common way to introduce fractional derivatives consists in first deriving a generalization to any order $\alpha>0$ of the integer-order integral, namely the fractional RL integral 
\begin{equation}\label{eq:RL_Integral}
 J_{t_0}^{\alpha} y(t) = \frac{1}{\Gamma(\alpha)} \int_{t_{0}}^{t} 
(t-\tau)^{\alpha-1} y(\tau) \du \tau 
	, \quad \forall t > t_0 ,
\end{equation}
and hence finding inverse operators which are therefore recognized as fractional derivatives. Actually, more than one operator performs the left-inversion of the integral (\ref{eq:RL_Integral}): one  is the \emph{RL fractional derivative}  
\begin{equation}\label{eq:RL_Derivative}
	\DR_{t_0}^{\alpha} y(t) := D^{m} J_{t_0}^{m-\alpha} y(t)  = \frac{1}{\Gamma(m-\alpha)} \frac{\du^{m}}{\du t^{m}} \int_{t_{0}}^{t} (t-\tau)^{m-\alpha-1} y(\tau) \du \tau,  \quad
	t > t_0 
\end{equation}
and a further one is the \emph{Caputo fractional derivative}
\begin{equation}\label{eq:Caputo_Derivative}
	\DC_{t_0}^{\alpha} y(t) := J_{t_0}^{m-\alpha} D^{m} y(t)  = \frac{1}{\Gamma(m-\alpha)} \int_{t_{0}}^{t} (t-\tau)^{m-\alpha-1} y^{(m)}(\tau) \du \tau , \quad
	t > t_0, 	
\end{equation}
where $m= \left\lceil \alpha \right \rceil$ is the smallest integer greater or equal to $\alpha$, $D^{m}$ and $y^{(m)}$ denote integer-order differentiation and $\Gamma(\beta)$ is the Euler-Gamma function. Under reasonable assumptions about the function $f$, one can indeed verify that $\DR_{t_0}^{\alpha} J_{t_0}^{\alpha} f = \DC_{t_0}^{\alpha} J_{t_0}^{\alpha} f = f$ (e.g., see \cite{Diethelm2010,KilbasSrivastavaTrujillo2006,SamkoKilbasMarichev1993}).

Although quite subtle, this procedure turns out to be convenient for devising a solid theory and, actually, the majority of the existing fractional derivatives are introduced by following an approach of this kind. However, obtaining fractional derivatives by inversion of the integral presents some drawbacks since it hides the contribution of the function to the left of the starting point $t_0$ which, as we will see, may have an important role in fractional calculus and for FDDEs in particular. It is therefore useful to present an alternative approach.

To this purpose, in the difference quotient definition of $n$-th integer-order derivatives 

\begin{equation}\label{eq:IntegerOrderDerivative}
	f^{(n)}(t) = \lim_{h\to0} \frac{1}{h^n} \sum_{j=0}^{n} \omega_{j}^{(n)} y(t-jh) 
	, \quad
	\omega_{j}^{(n)} = (-1)^j \binom{n}{j} 	,
\end{equation}
observe that coefficients $\omega_{j}^{(n)}$ can be formulated in terms of the Euler-Gamma function, as

\begin{equation}\label{eq:BinomialCoefficientsInteger}
	\binom{n}{j} =  \displaystyle\frac{n!}{j!(n-j)!} = \left\{\begin{array}{ll}
		\displaystyle\frac{\Gamma(n+1)}{j!\Gamma(n+1-j)} \quad & j=0,1,\dots,n , \\
		0 & j > n .\\
	\end{array}\right.
\end{equation}

To generalize (\ref{eq:IntegerOrderDerivative}) to any real order $\alpha>0$ it is necessary to first exploit the fact that $\omega_{j}^{(n)}=0$ for $j>n$, thus to be able to equivalently rewrite (\ref{eq:IntegerOrderDerivative}) as an infinite series, and hence replace the integer order $n$ with any real $\alpha>0$ to obtain

\begin{equation}\label{eq:FOD_GL}
	\DG^{\alpha} y(t) = \lim_{h\to0} \frac{1}{h^{\alpha}} \sum_{j=0}^{\infty} \omega_{j}^{(\alpha)} y(t-jh)
	, \quad 
	\omega_{j}^{(\alpha)} = \frac{(-1)^j\Gamma(\alpha+1)}{j! \Gamma(\alpha-j+1)} ,
\end{equation}
(see, for instance, \cite{GarrappaKaslikPopolizio2019} for a more complete discussion). The operator $\DG^{\alpha}$ is known as the Gr\"unwald-Letnikov (GL) fractional derivative, as it has been introduced almost simultaneously by Gr\"unwald \cite{Grunwald1867} and Letnikov \cite{Letnikov1868}, and represents one of the most straightforward ways to generalize the usual definition of the integer-order derivative.  It should be therefore quite natural to adopt (\ref{eq:FOD_GL}) as a standard way to introduce fractional derivatives but, unfortunately, some inconveniences arise:
\begin{enumerate}
\item the knowledge of $y(t)$ on the whole interval $[-\infty,t]$ is necessary to evaluate $\DG^{\alpha} y(t)$ and this may be a serious issue when the function describes a system state for which the history from $-\infty$ to the starting point is not available; 
\item the range of functions for which the series from (\ref{eq:FOD_GL}) converges, and hence $\DG^{\alpha} y(t)$ exists, is  restricted to bounded functions (when $0<\alpha<1$) or functions which do not increase too fast as $t\to -\infty$ \cite[\S 4.20]{SamkoKilbasMarichev1993}. 
\end{enumerate} 

For practical applications, some modifications are necessary to the operator $\DG^{\alpha}$ and the most frequently used approach is to select a starting point, which for convenience we choose at $0$, and force the function $y(t)$ to assume some suitably selected values in $(-\infty, 0)$. For instance, when $0<\alpha<1$, one may assume $y(t)=0$ or $y(t)=y_0$ for $t<0$ and, quite surprisingly, the fractional derivatives $\DR_{0}^{\alpha} y(t)$ or $\DC_{0}^{\alpha} y(t)$ are obtained. In particular, when $0<\alpha<1$ it is possible to show that
\[
	\DR_{0}^{\alpha} y(t) = \DG^{\alpha} \yR(t)
	, \quad	
	\yR(t) = \left\{ \begin{array}{ll}
		0 & t \in (-\infty, 0) \\
		y(t) \quad & t \ge 0 \
	\end{array}\right. 
\]
and
\begin{equation}\label{eq:ReplacmentGLCaputo}
	\DC_{0}^{\alpha} y(t) = \DG^{\alpha} \yC(t)
	, \quad	
	\yC(t) = \left\{ \begin{array}{ll}
		y_0 & t \in (-\infty, 0) \\
		y(t) \quad & t \ge 0 . \
	\end{array}\right.
\end{equation}

These replacements have the advantage of avoiding infinite series in the representation (\ref{eq:FOD_GL}) as well; this is obvious for the RL derivative and we observe that, since (e.g., see \cite{GarrappaKaslikPopolizio2019}) 
\begin{equation}\label{eq:BinomSeriesCoeffEquivalence}
	 \sum_{j=0}^{N} \omega_{j}^{(\alpha)} = - \sum_{j=N+1}^{\infty} \omega_{j}^{(\alpha)} ,
\end{equation}	
it is possible to obtain a more convenient representation for the Caputo derivative as well
\begin{equation}\label{eq:GL_Definition}
\DC_{0}^{\alpha} y(t) = \DG^{\alpha} \yC(t) = \lim_{h\to0} \frac{1}{h^{\alpha}} \sum_{j=0}^{\left\lfloor t/h \right\rfloor} \omega_{j}^{(\alpha)} \bigl[ y(t-jh) - y_0 \bigr] 
, \quad 0 < \alpha < 1 .
\end{equation}

The above discussion discloses the real nature of the RL and Caputo fractional derivatives: neither is actually a pure and straightforward generalization of the integer-order derivatives (\ref{eq:IntegerOrderDerivative}) but instead, they both stem from the the integer-order derivatives only after modifying the value of the function at the left of the selected starting point.

It is therefore mandatory to investigate the consequences of this peculiar nature of fractional-order derivatives on the solution of FDDEs in which the values for $-\tau \le t \le 0$ are determined by the initial condition $y(t)=\phi(t)$ as well.

To avoid inconsistencies, one could just impose initial functions compatible with the selected operator, for instance $\phi(t)=y_0$, $-\tau\le t\le 0$ for FDDEs with the Caputo derivative of order $0<\alpha<1$, but such a limitation of the initial data may be too restrictive for several models. Indeed, even if in many mathematical models which include time delays, the history is considered to be constant, non-constant history functions are also routinely encountered, such as in models arising from epidemiology and population dynamics \cite{kuang1993delay,smith2011introduction}. For instance, in epidemiological models, it is often assumed that the history of the number of infected individuals is a continuous, strictly increasing function $I_0(t)$ defined on the interval $[-\tau,0]$, such that $I_0(-\tau)=0$ and $I_0(0)>0$ \cite{hoppensteadt2006numerical}.

In the next sections we will therefore investigate two aspects. For a class of linear FDDEs of order $0<\alpha<1$ we will first derive in Section \ref{S:ExactSolution} the solution obtained with the operator $\DC_{0}^{\alpha}$, relying on the assumption $y(t)=y_0$ for $t\le 0$, and then, we will also derive the exact solution for different initial functions $\phi(t)$: in fact, we consider the case in which the initial function affects the solution but not the operator. Moreover, in Section \ref{S:GL_Modified} we will present the solution obtained by suitably modifying the fractional operator in order to comply with the initial data and force the solution to satisfy the initial condition $y(t)=\phi(t)$, $-\tau \le t \le 0$, and not just the condition $y(t)=y_0$, $t\le 0$, imposed by the Caputo derivative. 

We will not only show and compare the solutions obtained by the different approaches but we will also discuss the differences between the two operators and find the relationship connecting them.

\section{Exact solution of linear FDDEs}\label{S:ExactSolution}

With the aim of a better illustration of the effects of the initial function not just on the solution of the FDDE, but also on the nature of the fractional derivative, we derive here the exact solution for a family of linear FDDEs with the Caputo derivative (\ref{eq:Caputo_Derivative}) 
\begin{equation}\label{eq:FDDE_Linear}
\begin{cases}
	\DC^{\alpha}_0 y(t) = \lambda y(t-\tau) + f(t) &,~ t>0 \\
	y(t) = \phi(t)  &,~  -\tau < t \le 0 
\end{cases}
\end{equation}
where $\tau>0$, $\phi(t):[-\tau,0] \to \Rset$ and $f(t):[0,+\infty) \to \Rset$ are two given functions and $\lambda$ a constant (usually, but not necessarily, real) parameter. Since we restrict to FDDEs of order $0<\alpha<1$ with starting point $t_0=0$, the Caputo derivative is defined as $\DC_{0}^{\alpha} y(t) \coloneqq J^{\alpha}_{0} y'(t)$. 

Although this section is mainly intended to illustrate the effects of the initial function $\phi(t)$, we think that presenting exact solutions of linear FDDEs may have its own interest as well; analogous results were already presented in \cite{krol2011asymptotic,CermakHornicekKisela2016} for linear homogeneous FDDEs of non scalar type.

Based on Theorem \ref{thm.exp.bounded}, it is easy to see that the exponential boundedness of the solution of (\ref{eq:FDDE_Linear})  is guaranteed by the existence of two constants $M>0$ and $\beta>2\lambda$ such that 
\begin{equation}\label{eq:AssumptionLinearExpBound}
\int_0^t(t-r)^{\alpha-1}|f(r)| \du r \leq M E_\alpha(\beta t^\alpha), \quad\forall t\geq 0.
\end{equation} 

Therefore, if inequality \eqref{eq:AssumptionLinearExpBound} is fulfilled for the function $f(t)$, we can solve (\ref{eq:FDDE_Linear}) by means of the Laplace transform (LT). 

We preliminarily recall some function definitions and some results on the LT which will be useful in what follows. To this purpose we remember that for a continuous and exponential bounded function $f(t):[0,+\infty)\to\Rset$, its LT is
\[
	F(s) \coloneqq {\mathcal L} \bigl( f(t) \, ; \, s \bigr) = \int_0^{\infty} \eu^{-s t} f(t) \du t ,
\]
with $s$ belonging to the region in the complex plane where the above integral converges.

\subsection{Preliminary results}

For any $a\in \Rset$ the Heaviside unit step function is defined as
\[
	u_{a}(t) = \left\{\begin{array}{ll}
	1 \quad & t \ge  a \\
	0 \quad & t < a \
	\end{array}\right.
\]
and we can consider, for any real $\beta > -1$, its generalization
\[
	u_{a}^{[\beta]}(t) = \frac{(t-a)^\beta}{\Gamma(\beta+1)} u_{a}(t) = \left\{\begin{array}{ll}
		\displaystyle\frac{(t-a)^\beta}{\Gamma(\beta+1)} \quad & t \ge a \\
		0 \quad & t < a .\
		\end{array}\right.
\]
Clearly, it is $u_{a}(t) = u_{a}^{[0]}(t)$. We present here the following result on the LT of $u_{a}^{[\beta]}(t)$.

\begin{proposition}\label{prop:LTAntiderivativeStep}
Let $a>0$. For any $\beta>0$ the LT of $u_{a}^{[\beta]}(t)$ is 
\[
	{\mathcal L} \bigl( u_{a}^{[\beta]}(t) \, ; \, s \bigr) 
	= \frac{\eu^{-s a}}{s^{\beta+1}}
\]
\end{proposition}

\begin{proof}
It is immediate to observe that
\[
{\mathcal L} \bigl( u_{a}^{[\beta]}(t) \, ; \, s \bigr) = \int_0^{\infty} \eu^{-s t} \frac{(t-a)^{\beta}}{\Gamma(\beta+1)} u_{a}(t) \du t
= \int_{a}^{\infty}\eu^{-s t} \frac{(t-a)^{\beta}}{\Gamma(\beta+1)} \du t
\]
and after the change of variable $v=t-a$ it is
\[
	{\mathcal L} \bigl( u_{a}^{[\beta]}(t) \, ; \, s \bigr) 
	= \eu^{-s a} \int_{0}^{\infty} \eu^{-s v} \frac{v^\beta}{\Gamma(\beta+1)} \du v = \eu^{-s a} {\mathcal L} \Bigl( \frac{v^\beta}{\Gamma(\beta+1)} \, ; \, s \Bigr)  = \frac{\eu^{-s a}}{s^{\beta+1}}
\]
which concludes the proof.
\end{proof}

For the sake of completeness, we recall here the following results (e.g., see \cite{CermakHornicekKisela2016,KaslikSivasundaram2012}).

\begin{proposition}\label{prop:LT_DDE_Integer}
Let $\tau>0$, $y(t):[-\tau,+\infty) \to \Rset$ and $Y(s)$ its LT. Then
\begin{itemize}
\item ${\mathcal L} \bigl( y(t-\tau)u_{\tau}(t) \, ; \, s \bigr) = \eu^{-s \tau} Y(s)$;
\item ${\mathcal L} \bigl( y(t-\tau) \, ; \, s \bigr) = \eu^{-s \tau} \hat{Y}_{\tau}(s) + \eu^{-s \tau} Y(s), \quad
\hat{Y}_{\tau}(s) = \int_{-\tau}^{0} \eu^{-st } y(t) \du t$. 
\end{itemize}
\end{proposition}

It is well-known that the LT of the first-order derivative of $y(t)$ is ${\mathcal L} \bigl( y'(t) \, ; \, s \bigr) = s Y(s) - y_0$, with $y_0=y(0)$. To obtain the LT of the Caputo derivative $\DC_0^{\alpha} y(t)$ it is possible to start from the RL integral (\ref{eq:RL_Integral}), which is actually the convolution of the two functions $t^{\alpha-1}/\Gamma(\alpha)$ and $y(t)$, and its LT is easily evaluated as ${\mathcal L} \bigl(J^{\alpha}_0 y(t) \, ; \, s \bigr) = Y(s)/s^{\alpha}$. Therefore, the LT of $\DC_0^{\alpha} y(t) \coloneqq J^{1-\alpha} y'(t)$ follows standard rules in LT calculus and hence 
\begin{equation}\label{eq:LT_Caputo}
{\mathcal L} \bigl(\DC_0^{\alpha} y(t) \, ; \, s \bigr) = \frac{1}{s^{1-\alpha}} {\mathcal L} \bigl(y'(t) \, ; \, s \bigr)
= \frac{ s Y(s) - y_0}{s^{1-\alpha}} = s^{\alpha} Y(s) - s^{\alpha-1} y_0 .
\end{equation}

\begin{remark}
Although, as discussed in Section \ref{S:GeneralizingDerivatives}, the Caputo derivative of order $0<\alpha<1$ is equivalent to the operator obtained from the GL derivative when the function is forced to assume constant value $y_0$ for $t<0$, the LT of $\DC_0^{\alpha} y(t)$ is obtained only on the basis of the integral representation and is therefore independent of any further change imposed on the behavior of the function as consequence of the initial condition in (\ref{prop:LT_DDE_Integer}). Obviously, the same invariance is not expected when the derivative is obtained from the approach outlined by (\ref{eq:ReplacmentGLCaputo}) which instead reflects changes of the function to the left of the origin. Thus, the solution we are going to compute will present a certain level of contradiction: the operator considers the solution $y(t)=y_0$ for $t<0$ but at the same time, the equation imposes $y(t)=\phi(t)$, $t\in[-\tau,0]$, with the initial function $\phi(t)$ which could differ from $y_0$.
\end{remark}

\subsection{Exact solution of linear fractional-order DDEs}

Just for notational convenience, we introduce the following generalized integral 
\[
	\mathcal{J}_{0,\tau,\lambda}f(t) \coloneqq \sum_{k=0}^{\left\lfloor t/\tau \right\rfloor} \frac{\lambda^{k} }{\Gamma\bigl(\alpha k + \alpha \bigr)} \int_{0}^{t-k\tau} (t-k \tau - r )^{ \alpha k+ \alpha-1}  f(r) \du r.
\]
and we provide a first general representation of the exact solution of (\ref{eq:FDDE_Linear}).

\begin{proposition}\label{prop:DDE_Frac_ExSol}
Let $\tau > 0$. For any $t\ge 0$ the exact solution of the linear FDDE (\ref{eq:FDDE_Linear}) can be expressed in the form
\[
	y(t) = y_0 + \sum_{k=0}^{\left\lfloor t/\tau \right\rfloor} \frac{\lambda^{k+1}}{\Gamma((k+1)\alpha+1)} \int_{0}^{t} (t- k \tau  -r) ^{(k+1) \alpha} p_{\tau}(r;\phi) \du r   + \mathcal{J}_{0,\tau,\lambda}f(t) ,
\]
where $y_0 = \phi(0)$ and the function $p_{\tau}(t;\phi)$ is the inverse LT of $\eu^{-s\tau} P_{\tau} (s;\phi)$, with
\[
P_{\tau} (s;\phi) = \phi(0) + s \int_{-\tau}^{0} \eu^{-s r} \phi(r)  \du r .
\]
\end{proposition}

\begin{proof}
By exploiting the results from Proposition \ref{prop:LT_DDE_Integer} and the LT (\ref{eq:LT_Caputo}) of the Caputo's derivative, it is possible to represent the LT of the solution of (\ref{eq:FDDE_Linear}) as 
\[
	s^{\alpha}Y(s) - s^{\alpha-1} y_0 = \lambda \eu^{-s \tau} \hat{\Phi}_{\tau}(s) + \lambda \eu^{-s \tau} Y(s) + F(s)
\]
with $Y(s)$ and $F(s)$ the LT of $y(t)$ and $f(t)$ respectively. It is hence immediate to show that
\[
	\begin{aligned}
	Y(s) 
	&= \left( \frac{1}{ s^{\alpha} - \lambda \eu^{-s \tau}} \right) \left[ s^{\alpha-1} y_0 + \lambda \eu^{-s \tau}\int_{-\tau}^{0} \eu^{-sr } \phi(r) \du r + F(s) \right] \\
	&=\frac{1}{s^{\alpha}} \left(1 - \lambda \frac{\eu^{-s \tau}}{s^{\alpha}} \right)^{-1}  \left[ s^{\alpha-1}\left( 1 - \lambda \frac{\eu^{-s\tau}}{s^{\alpha}}\right) y_0 + \lambda \frac{\eu^{-s\tau}}{s} P_{\tau} (s;\phi)  \right]  + \frac{1}{s^{\alpha}} \left(1 - \lambda \frac{\eu^{-s \tau}}{s^{\alpha}} \right)^{-1} F(s) \\
	&= \frac{1}{s}y_0 + \left(1 - \lambda \frac{\eu^{-s \tau}}{s^{\alpha}} \right)^{-1} \lambda \frac{\eu^{-s\tau}}{s^{1+\alpha}} P_{\tau} (s;\phi)   + \frac{1}{s^{\alpha}} \left(1 - \lambda \frac{\eu^{-s \tau}}{s^{\alpha}} \right)^{-1} F(s) . \\
	\end{aligned}
\]

For sufficiently large $|s|$ consider the series expansion
\[
	\left(1 - \lambda \frac{\eu^{-s \tau}}{s^{\alpha}} \right)^{-1} = \sum_{k=0}^{\infty} \lambda^k \frac{\eu^{-s\tau k}}{s^{k\alpha}}
\]
in order to represent the LT of the solution of (\ref{eq:FDDE_Linear}) as
\begin{equation}\label{eq:FDDE_LT_Sol}
	Y(s) 
	= \frac{1}{s}y_0 + \sum_{k=0}^{\infty} \lambda^{k+1} \frac{\eu^{-s\tau k}}{s^{(k+1)\alpha+1}} \eu^{-s\tau} P_{\tau} (s;\phi)    + \sum_{k=0}^{\infty} \lambda^k \frac{\eu^{-s\tau k}}{s^{(k+1)\alpha}} F(s) 
\end{equation}
and, after inverting the LT and using the results in Proposition \ref{prop:LTAntiderivativeStep}, we obtain
\[
	y(t) = y_0 + \sum_{k=0}^{\infty} \lambda^{k+1} \int_{0}^{t} u_{k\tau}^{[(k+1)\alpha]}(t-r) p_{\tau}(r;\phi) \du r   + \sum_{k=0}^{\infty} \lambda^k \int_{0}^{t} u_{k\tau}^{[(k+1)\alpha-1]}(t-r) f(r) \du r .
\]

The conclusion now follows by replacing each function $u_{a}^{[\beta]}$ with the corresponding power function, and in the proper domain in which it does not vanish, thus to be able to truncate the infinite series.
\end{proof}

To better observe the influence of the initial function on the solution and provide easily evaluable formulas, we must select the function $\phi(t)$. For a free value $y_0 \in \Rset$, throughout this paper we will consider two exemplifying functions:
\begin{enumerate}
	\item the constant function $\phi(t)\equiv y_0$ which agrees with the assumption made on the solution $y(t)$ to obtain the equivalence between the Caputo and the GL derivatives;
	\item a first order polynomial $\phi(t)=(t/\tau + 1)y_0$ which instead imposes a different behavior to the solution $y(t)$, compared to the one which ensures the equivalence between the Caputo and the GL derivatives.
\end{enumerate}

\begin{corollary}\label{cor:FDDE_Lin_Ex_FirstInt_Case1}
Let $\tau > 0$ and $\phi(t)\equiv y_0$, $t \in [-\tau,0]$, for some $y_0 \in \Rset$. For any $t \ge 0$ the exact solution of the linear FDDE (\ref{eq:FDDE_Linear}) is
\begin{equation}\label{eq:FDDE_Lin_Ex_FirstInt_Case1}
	y(t) =  \sum_{k=0}^{\left\lfloor t/\tau \right\rfloor+1} \frac{\lambda^{k} (t+\tau-k \tau)^{\alpha k}}{\Gamma\bigl(\alpha k+1\bigr)}  y_0 
	+ \mathcal{J}_{0,\tau,\lambda}f(t) .
\end{equation}
\end{corollary}
\begin{proof}  
It is immediate in this case to verify that 
\[
	\eu^{-s\tau}  P_{\tau} (s;\phi)
	= \left[ \eu^{-s\tau} + \eu^{-s\tau} s \int_{-\tau}^{0} \eu^{-s r}  \du r \right] y_0
	= \left[\eu^{-s\tau} + \eu^{-s\tau} \left(  \eu^{s \tau} - 1 \right)  \right] y_0 = y_0
\]
and therefore the representation (\ref{eq:FDDE_LT_Sol}) of the LT of the solution of (\ref{eq:FDDE_Linear}) is 
\[
	Y(s) 
	= \left[ \frac{1}{s}y_0 + \sum_{k=0}^{\infty} \lambda^{k+1} \frac{\eu^{-s\tau k}}{s^{(k+1)\alpha+1}} y_0  \right]  + \sum_{k=0}^{\infty} \lambda^k \frac{\eu^{-s\tau k}}{s^{(k+1)\alpha}} F(s) .
\]
Thanks to Proposition \ref{prop:LTAntiderivativeStep}, the inversion of the LT leads to
\[
	y(t) =  \left[ 1 +  \sum_{k=0}^{\infty} \lambda^{k+1} u_{\tau k}^{[(k+1)\alpha]}(t) \right] y_0 + \sum_{k=0}^{\infty} \lambda^k \int_{0}^{t} u_{k\tau}^{[(k+1)\alpha-1]}(t-r) f(r) \du r
\]
and the corollary follows by replacing the functions $u_{a}^{[\beta]}(t)$ with the corresponding powers in the proper domain and after reorganizing some terms.
\end{proof}

\begin{figure}[tph]
\centering
\includegraphics[width=0.95\textwidth]{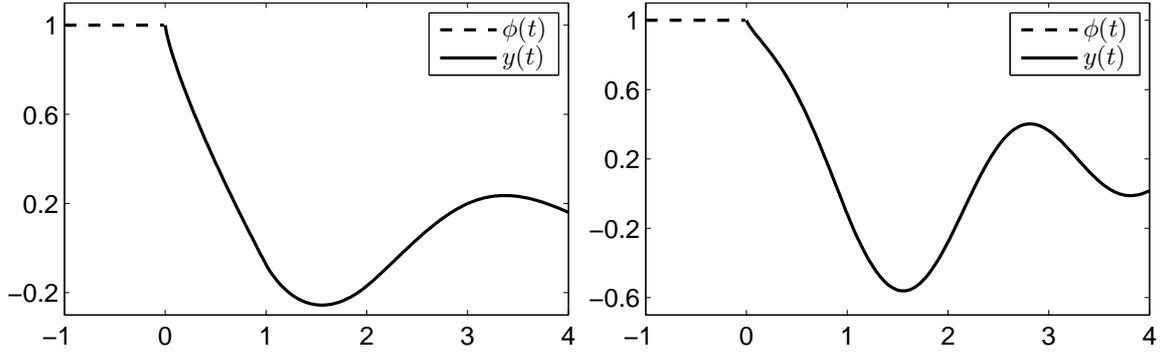}
\caption{Plot of the solution (\ref{eq:FDDE_Lin_Ex_FirstInt_Case1}) for $\lambda=-1.0$, $\tau=1.0$,  $\phi(t)=y_0$ and $f(t)=0$ (left plot) and $f(t) = \frac{1}{2} \cos 3 t$ (right plot)\label{Fig_FDDE_ExactSolution_Linear_1}}
\end{figure}

In Figure \ref{Fig_FDDE_ExactSolution_Linear_1} we show the solution (\ref{eq:FDDE_Lin_Ex_FirstInt_Case1}) for the homogeneous case $f(t)=0$ (left plot) and for the forcing function $f(t) = \frac{1}{2} \cos 3 t$ (right plot); the values $\alpha=0.8$, $\lambda=-1$, $\tau=1$ and $y_0=1$ have been used. The integral $\mathcal{J}_{0,\tau,\lambda}f(t)$ is evaluated by exploiting the exact formula of the RL integral of the cosine function in terms of the Mittag-Leffler (ML) function (e.g., see \cite{GarrappaKaslikPopolizio2019}), with the ML function evaluated by means of the Matlab code developed in \cite{Garrappa2015}. 

\begin{remark}
Observe that (\ref{eq:FDDE_Lin_Ex_FirstInt_Case1}) straightforwardly generalizes the well-known variation-of-constants formula (see, e.g., \cite{GarrappaPopolizio2011_CAMWA})
\[
	y(t) = e_{\alpha,1}\bigl(t;\lambda\bigr) y_0 + \int_{0}^{t} e_{\alpha,\alpha}\bigl((t-r);\lambda\bigr) f(r) \du r
\]
for the linear (non-delayed) FDE $\DC_{0}^{\alpha} y(t) = \lambda y(t) + f(t)$, where $e_{\alpha,\beta}(t;\lambda)$ is the generalized ML function
\[
	e_{\alpha,\beta}(t;\lambda) = t^{\beta-1} E_{\alpha,\beta}(t^{\alpha} \lambda)
	, \quad
	E_{\alpha,\beta}(z) = \sum_{k=0}^{\infty}\frac{z^k}{\Gamma(\alpha k + \beta)} .
\]
\end{remark}

\begin{corollary}\label{cor:FDDE_Lin_Ex_FirstInt_Case2}
Let $\tau > 0$ and $\phi(t) = \bigl(t/\tau + 1)y_0$, $t \in [-\tau,0]$, for some $y_0 \in \Rset$. For any $t \ge 0$ the exact solution of the linear FDDE (\ref{eq:FDDE_Linear}) is
\begin{equation}\label{eq:FDDE_Lin_Ex_FirstInt_Case2}
	\begin{aligned}
	y(t) = \frac{1}{\tau}  \left[ \sum_{k=0}^{\left\lfloor t/\tau \right\rfloor+1} \frac{\lambda^{k} (t+\tau-\tau k)^{\alpha k+1}}{\Gamma\bigl(\alpha k+2\bigr)} - \sum_{k=0}^{\left\lfloor t/\tau \right\rfloor} \frac{\lambda^{k} (t-\tau k)^{\alpha k+1}}{\Gamma\bigl(\alpha k+2\bigr)}  \right] y_0 +
\mathcal{J}_{0,\tau,\lambda}f(t)	\end{aligned}
\end{equation}
\end{corollary}
\begin{proof}
By standard integration we preliminarily observe that when $\phi(t)= \bigl(t/\tau + 1)y_0$ it is
\[
\begin{aligned}
	\eu^{-s\tau} & P_{\tau} (s;\phi)
	= e^{-s\tau}\left[ 1 + \frac{1}{\tau} \int_{-\tau}^{0}  r \eu^{-s r} s \du r +  \int_{-\tau}^{0} \eu^{-s r} s\du r \right] y_0\\
	&= e^{-s\tau}\left[ 1 + \frac{1}{s\tau} \left( -s\tau e^{s\tau} +e^{s\tau}-1 \right)  +  \eu^{s\tau} - 1  \right] y_0 	
	= \left[-\frac{\eu^{-s\tau}}{s \tau} + \frac{1}{s\tau} \right] y_0
\end{aligned}
\]
and hence the representation (\ref{eq:FDDE_LT_Sol}) of the LT of the solution of (\ref{eq:FDDE_Linear}) becomes
\[
	\begin{aligned}
	Y(s) 
	&= \left[ \frac{1}{s}y_0 + \sum_{k=0}^{\infty} \lambda^{k+1} \frac{\eu^{-s\tau k}}{s^{(k+1)\alpha+1}} \left[-\frac{\eu^{-s\tau}}{s \tau} + \frac{1}{s\tau} \right] y_0  \right]  + \sum_{k=0}^{\infty} \lambda^k \frac{\eu^{-s\tau k}}{s^{(k+1)\alpha}} F(s) \\
	&= \left[ \frac{1}{s} - \frac{1}{\tau} \sum_{k=0}^{\infty} \lambda^{k+1} \frac{\eu^{-s\tau(k+1)}}{s^{(k+1)\alpha+2}} + \frac{1}{\tau}\sum_{k=0}^{\infty} \lambda^{k+1} \frac{\eu^{-s\tau k}}{s^{(k+1)\alpha+2}} \right] y_0 + \sum_{k=0}^{\infty} \lambda^k \frac{\eu^{-s\tau k}}{s^{(k+1)\alpha}} F(s) . \\
	\end{aligned}
\]
By inversion of the LT it is therefore 
\begin{equation}\label{eq:FDDE_Lin_Ex_FirstInt_Case2_StepFunction}
	\begin{aligned}
	y(t) &= \left[ 1 - \frac{1}{\tau} \sum_{k=0}^{\infty} \lambda^{k+1} u_{\tau k+\tau}^{[(k+1)\alpha+1]}(t) + \frac{1}{\tau} \sum_{k=0}^{\infty} \lambda^{k+1} u_{\tau k}^{[(k+1)\alpha+1]}(t) \right] y_0 \\
	&\quad + \sum_{k=0}^{\infty} \lambda^k \int_{0}^{t} u_{k\tau}^{[(k+1)\alpha-1]}(t-r) f(r) \du r \\
	\end{aligned}
\end{equation}
and, again, a proper replacement of the functions $u_{a}^{[\beta]}(t)$  with the corresponding power functions, together with the reorganization of some of the summations, allows to conclude the proof. 
\end{proof}

\begin{figure}[tph]
\centering
\includegraphics[width=0.95\textwidth]{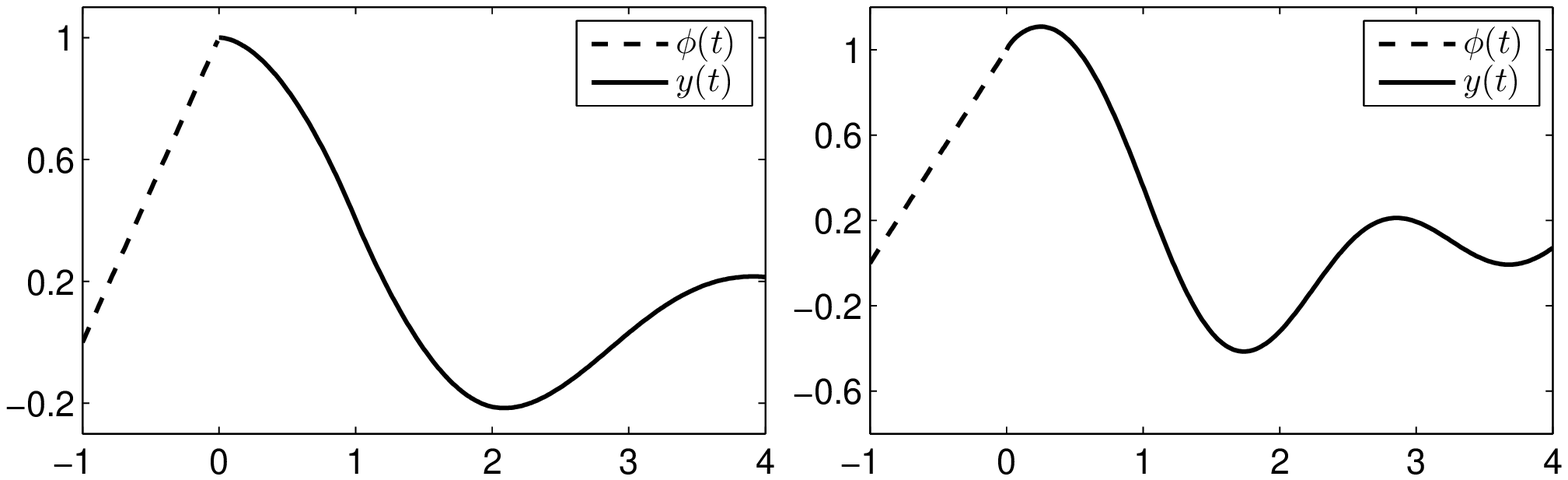}
\caption{Plot of the solution (\ref{eq:FDDE_Lin_Ex_FirstInt_Case2}) for $\lambda=-1.0$, $\tau=1.0$,  $\phi(t)=\bigl(t/\tau + 1)y_0$ and $f(t)=0$ (left plot) and $f(t) = \frac{1}{2} \cos 3 t$ (right plot)\label{Fig_FDDE_ExactSolution_Linear_2}}
\end{figure}

The solution (\ref{eq:FDDE_Lin_Ex_FirstInt_Case2}) is presented in Figure \ref{Fig_FDDE_ExactSolution_Linear_2} for the same data $\lambda$, $\tau$ and $f(t)$ used to plot the solution (\ref{eq:FDDE_Lin_Ex_FirstInt_Case1}) in Figure \ref{Fig_FDDE_ExactSolution_Linear_1} (the selected data are anyway reported in the caption).

\section{Incorporating the initial function in the GL fractional derivative}\label{S:GL_Modified}

Let us consider now the FDDE (\ref{eq:FDDE_General}) with a fractional derivative obtained after generalizing the limit of the difference quotient defining the integer-order derivative, according to the process already described in Section \ref{S:GeneralizingDerivatives}. 

Because of the initial condition $y(t)=\phi(t)$, $-\tau\le t \le0$, to avoid the infinite memory and ensure the convergence, we cannot simply force $y(t)=y_0$ for $t \le 0$ as in (\ref{eq:ReplacmentGLCaputo}). We can consider the following three different options which resemble the same process leading to the Caputo derivative and, at the same time, fulfill the initial condition: 


\begin{enumerate}
	\item \emph{assuming $y(t)=\phi(t)$ for any $t \le 0$}: this option may apply when $\phi(t)$ is defined on $(-\infty,0]$ and not only on $[-\tau,0]$ but, unfortunately, it is in general an unfeasible option since for several functions $\phi(t)$ the series from (\ref{eq:FOD_GL}) might not converge;
	\item \emph{assuming $y(t)=\phi(t)$ for $-\tau \le t \le 0$ and $y(t)=y_0$ for $t < -\tau$}: although this choice would appear as a natural and light adjustment of the Caputo derivative for this problem, it would introduce an undesirable discontinuity when $\phi(0)\not=\phi(-\tau)$ (we must remember that the Caputo derivative has been introduced just to regularize a similar discontinuity induced by the RL derivative);
	\item \emph{assuming $y(t)=\phi(t)$ for $-\tau \le t \le 0$ and $y(t)=\phi(-\tau)$ for $t < -\tau$}: in our opinion, this last option appears to be the most reasonable, since it ensures the convergence of the series from (\ref{eq:FOD_GL}) and the fulfillment of the initial condition without introducing an unwanted discontinuity.
\end{enumerate}

To apply the above third option we just replicate the process for the construction of the Caputo derivative (i.e. replacing in the GL derivative the function with another function), but instead of (\ref{eq:ReplacmentGLCaputo}) we operate the replacement
\[
	\Dphi_{{\scriptscriptstyle 0}}^{\alpha} y(t) \coloneqq \DG^{\alpha} \yphi(t)
	, \quad	
	\yphi(t) = \left\{ \begin{array}{ll}
		\phitau(t ) & t \leq 0 \\
		y(t) \quad & t > 0 \
	\end{array}\right.
\]
where
\[
\phitau(t ) = \phi(-\tau) + \bigl[ \phi(t) - \phi(-\tau) \bigr]u_{-\tau}(t) = \left\{ \begin{array}{ll}
		\phi(-\tau) & t \in (-\infty, -\tau) \\
		\phi(t) & t \in [-\tau,0] \\
	\end{array}\right. .
\]

The resulting operator is therefore given by 
\begin{equation}\label{eq:GL_Deriv_DDE}
	\begin{aligned}
	\Dphi_{{\scriptscriptstyle 0}}^{\alpha} y(t)  
	&= \lim_{h\to0} \frac{1}{h^{\alpha}} \sum_{j=0}^{\infty} \omega_{j}^{(\alpha)} \yphi(t)(t-jh) \\
	&= \lim_{h\to0} \frac{1}{h^{\alpha}} \left[ \sum_{j=0}^{J_1} \omega_{j}^{(\alpha)} y(t-jh) + \sum_{j=J_1+1}^{J_2} \omega_{j}^{(\alpha)} \phi(t-jh) + \sum_{j=J_2+1}^{\infty} \omega_{j}^{(\alpha)} \phi(-\tau) \right] \\
	\end{aligned}
\end{equation}
where	
\[
	J_1 \coloneqq J_1(t,h) = \left\lfloor \frac{t}{h} \right\rfloor
	, \quad
	J_2 \coloneqq J_2(t,h,\tau) = \left\lfloor \frac{t+\tau}{h} \right\rfloor .
\]

The infinite series can be avoided by applying (\ref{eq:BinomSeriesCoeffEquivalence}) for $N=J_2$. We therefore obtain the more convenient representation of $\Dphi_{{\scriptscriptstyle 0}}^{\alpha} y(t)$ 

\begin{equation}\label{eq:FracDerviPhi_CalcualtingFormula}
	\Dphi_{{\scriptscriptstyle 0}}^{\alpha} y(t)  
	= \lim_{h\to0} \frac{1}{h^{\alpha}} \left[ \sum_{j=0}^{J_1} \omega_{j}^{(\alpha)} \bigl( y(t-jh) - \phi(-\tau) \bigr) + \! \! \sum_{j=J_1+1}^{J_2} \! \! \omega_{j}^{(\alpha)} \bigl( \phi(t-jh) - \phi(-\tau) \bigr) \right] .
\end{equation}

Due to the change operated in the replacement of $y(t)$ in $[-\tau,0]$ we can no longer expect the equivalence between $\Dphi_{{\scriptscriptstyle 0}}^{\alpha}$ and $\DC_{0}^{\alpha}$. It is useful to explore the relationship between these two operators in order to disclose the true nature of the new operator $\Dphi_{{\scriptscriptstyle 0}}^{\alpha}$.

\begin{theorem}\label{eq:EquivalenceNewOperator} 
Let $\tau>0$. Then
\[
	\Dphi_{{\scriptscriptstyle 0}}^{\alpha} y(t) = \DC_{0}^{\alpha} y(t) + \DC_{-\tau}^{\alpha} \phi(t)  - \DC_{0}^{\alpha}  \phi(t) .
\]
\end{theorem}

\begin{proof}
It is sufficient to add and subtract the same terms in (\ref{eq:FracDerviPhi_CalcualtingFormula}) and reorganize the resulting summations to observe that
\[
	\begin{aligned}
	\Dphi_{{\scriptscriptstyle 0}}^{\alpha} y(t)  
	&= \lim_{h\to0} \frac{1}{h^{\alpha}} \Biggl[ \sum_{j=0}^{J_2} \omega_{j}^{(\alpha)} \bigl( y(t-jh) - \phi(-\tau) \bigr) + \sum_{j=0}^{J_2} \omega_{j}^{(\alpha)} \bigl( \phi(t-jh) - y(t-jh) \bigr)  - \\
	&\phantom{\lim_{h\to0} \frac{1}{h^{\alpha}} \Biggl[}
	  \quad \quad - \sum_{j=0}^{J_1} \omega_{j}^{(\alpha)} \bigl( \phi(t-jh) - y(t-jh) \bigr) \Biggr] . \\
	\end{aligned}
\]

The first summation is clearly the Caputo fractional derivative of $y(t)$ with initial point at $-\tau$. Moreover, because of the initial conditions, it is $\phi(-\tau) - y(-\tau) = 0$ and $\phi(0) - y(0) = 0$ and also the second and third summations are Caputo derivatives of $\phi(t) - y(t)$ with starting point respectively at $-\tau$ and $0$. Therefore it is 
\[
	\Dphi_{{\scriptscriptstyle 0}}^{\alpha} y(t) 
	= \DC_{-\tau}^{\alpha} y(t) + \DC_{-\tau}^{\alpha} \bigl( \phi(t) - y(t) \bigr) - \DC_{0}^{\alpha} \bigl( \phi(t) - y(t) \bigr)  
\]
and the proof follows from the linearity of $\DC_{-\tau}^{\alpha}$ and $\DC_{0}^{\alpha}$.
\end{proof}



From Theorem \ref{eq:EquivalenceNewOperator} we infer that solving the FDDE (\ref{eq:FDDE_General}) with the fractional derivative $\Dphi_{{\scriptscriptstyle 0}}^{\alpha}$ obtained as a generalization of  the integer-order derivative which, at the same time, satisfies the initial condition of the FDDE, is equivalent to solving the equation 
\begin{equation}\label{eq:FDDE_General_Modified}
	\begin{cases}
	\DC^{\alpha}_{0} y(t) = g(t,y(t),y(t-\tau))  + \DC_{0}^{\alpha}  \phi(t) - \DC_{-\tau}^{\alpha} \phi(t) &,~t>0 \\
	y(t) = \phi(t) &, ~ -\tau \le t \le 0 .\\
	\end{cases}
\end{equation}

For the linear FDDE (\ref{eq:FDDE_Linear}) we are now able to compare the exact solution when using the derivative $\DC_{0}^{\alpha}$ with the exact solution obtained by emplying the modified derivative $\Dphi_{{\scriptscriptstyle 0}}^{\alpha}$. It is indeed sufficient to add the \emph{corrective term} $\DC_{0}^{\alpha}  \phi(t) - \DC_{-\tau}^{\alpha} \phi(t)$ to the forcing term $f(t)$ and hence use the results presented in Section \ref{S:ExactSolution}. We perform this comparison for the two examples of initial data $\phi(t)$ we are considering throughout this paper.



\subsection{Constant initial function}

Whenever $\phi(t)\equiv y_0$ the corrective term $\DC_{0}^{\alpha}  \phi(t) - \DC_{-\tau}^{\alpha} \phi(t)$ clearly vanishes and hence $\Dphi_{{\scriptscriptstyle 0}}^{\alpha} y(t)  = \DC_{0}^{\alpha} y(t)$. 

Assuming a constant history of the solution before the initial time is actually the assumption made in the case of the Caputo derivative in (\ref{eq:ReplacmentGLCaputo}) and the two operators coincide.

\subsection{First degree polynomial initial function}

When $\phi(t) = (t/\tau +1)y_0$ we can easily evaluate that the corrective term as
\begin{equation}\label{eq:CorrectiveTerm_Case2}
	\DC_{0}^{\alpha}  \phi(t) - \DC_{-\tau}^{\alpha} \phi(t)  = \frac{y_0}{\tau\Gamma(2-\alpha)} \Bigl(  t^{1-\alpha}-(t+\tau)^{1-\alpha}  \Bigr)
\end{equation}
and the action of $\Dphi_{{\scriptscriptstyle 0}}^{\alpha}$ is expected to differ from that of $\DC_{0}^{\alpha}$.

To evaluate the exact solution of (\ref{eq:FDDE_Linear}) under the operator $\Dphi_{{\scriptscriptstyle 0}}^{\alpha}$ it is sufficient to consider the solution under the operator $\DC_{0}^{\alpha}$ after adding the corrective term (\ref{eq:CorrectiveTerm_Case2}) to the source term $f(t)$. To this purpose we first introduce the following preliminary results.

\begin{lemma}\label{eq:RL_Integral_PowerFunctions}
Let $\beta > -1$ and $\alpha > 0$. Then
\begin{enumerate}
	\item $\displaystyle\frac{1}{\Gamma(\beta)}\int_{0}^{t}(t-r)^{\beta-1} \frac{r^{1-\alpha}}{\Gamma(2-\alpha)} \du r = \frac{t^{1-\alpha+\beta}}{\Gamma(2-\alpha+\beta)}$ ;
	\item $\displaystyle \frac{1}{\Gamma(\beta)}\int_{0}^{t}(t-r)^{\beta-1} \frac{(r+\tau)^{1-\alpha}}{\Gamma(2-\alpha)} \du r = \frac{(t+\tau)^{1-\alpha+\beta}}{\Gamma(2-\alpha+\beta)} -  \frac{(t+\tau)^{1-\alpha+\beta}}{\Gamma(2-\alpha+\beta)} I_{\frac{\tau}{t+\tau}}(2-\alpha,\beta) $,
\end{enumerate}
where $I_{x}\bigl(a,b\bigr)$ is the incomplete beta function \cite[Formula 6.6.2]{AbramowitzStegun1964}
\[
	I_{x}\bigl(a,b\bigr) = \frac{1}{B(a,b)}\int_0^x r^{a-1}(1-r)^{b-1} \du r .
\]

\end{lemma}

\begin{proof}
The first point is an immediate consequence of the formula for the RL integral of power functions (e.g., see \cite[Example 2.1]{Diethelm2010}). For the second point,  simple changes of variables allow to observe that 
\[
	\begin{aligned}
	\frac{1}{\Gamma(\beta)}\int_{0}^{t}&(t-r)^{\beta-1} \frac{(r+\tau)^{1-\alpha}}{\Gamma(2-\alpha)} \du r 
	= \frac{1}{\Gamma(\beta)}\int_{\tau}^{t+\tau}(t+\tau-r)^{\beta-1} \frac{r^{1-\alpha}}{\Gamma(2-\alpha)} \du r =   \\ 
	&= \frac{1}{\Gamma(\beta)}\int_{0}^{t+\tau}(t+\tau-r)^{\beta-1} \frac{r^{1-\alpha}}{\Gamma(2-\alpha)} \du r - \frac{1}{\Gamma(\beta)}\int_{0}^{\tau}(t+\tau-r)^{\beta-1} \frac{r^{1-\alpha}}{\Gamma(2-\alpha)} \du r =   \\ 
		&= \frac{1}{\Gamma(\beta)}\int_{0}^{t+\tau}(t+\tau-r)^{\beta-1} \frac{r^{1-\alpha}}{\Gamma(2-\alpha)} \du r - \frac{(t+\tau)^{1-\alpha+\beta}}{\Gamma(\beta)\Gamma(2-\alpha)}\int_{0}^{\frac{\tau}{t+\tau}}(1-r)^{\beta-1} r^{1-\alpha} \du r    \\ 
	\end{aligned}
\]
and the proof follows after applying the first point and the definition of the incomplete beta function.
\end{proof}

From the above Lemma we first observe that when we add the corrective term $\DC_{0}^{\alpha}  \phi(t) - \DC_{-\tau}^{\alpha} \phi(t)$ to the source term of the linear FDDE (\ref{eq:FDDE_Linear}), since  
\[
	\begin{aligned}
	\int_0^t(t-r)^{\alpha-1}&\bigl|\DC_{0}^{\alpha}  \phi(r) - \DC_{-\tau}^{\alpha} \phi(r)  \bigr| \du r \\
	&=|y_0|\Gamma(\alpha)\left[1-\frac{t+\tau}{\tau}I_{\frac{\tau}{t+\tau}}(2-\alpha,\alpha)\right]
	\leq |y_0|\Gamma(\alpha),\quad\forall t\geq 0 \\
	\end{aligned}
\]
it is sufficient to chose  $M=|y_0|\Gamma(\alpha)$, as $E_{\alpha}(\beta t^\alpha)\geq 1$, for any $\beta>0$ and $t\geq 0$, in order to ensure that the assumption (\ref{eq:AssumptionLinearExpBound}) is verified.

Let us now denote by $\hat{y}(t)$ the solution of the linear FDDE (\ref{eq:FDDE_Linear}) with the Caputo derivative $\DC^{\alpha}_0$ replaced by the fractional derivative $\Dphi^{\alpha}_0$, namely 
\begin{equation}\label{eq:FDDE_Linear_ModifiedOperator}
\begin{cases} 
	\Dphi^{\alpha}_0 \hat{y}(t) = \lambda \hat{y}(t-\tau) + f(t) &,~t>0, \\
	\hat{y}(t) = \phi(t) &,~  -\tau < t \le 0 ,
\end{cases}
\end{equation}
which in turn, in view of (\ref{eq:FDDE_General_Modified}), is solution of the equivalent FDDE
\begin{equation}\label{eq:FDDE_Linear_ModifiedOperator_Equivalent}
	\begin{cases}
	\DC^{\alpha}_0 \hat{y}(t) = \lambda \hat{y}(t-\tau) + f(t) + \DC_{0}^{\alpha}  \phi(t) - \DC_{-\tau}^{\alpha} \phi(t) &,~t>0,\\
	\hat{y}(t) = \phi(t) , &,~ -\tau < t \le 0. 
\end{cases}
\end{equation}

It is possible to provide an explicit representation of the difference between $\hat{y}(t)$ and the solution $y(t)$ of (\ref{eq:FDDE_Linear}) by means of the following result.

\begin{proposition}\label{prop:DifferenceExactSolutions}
Let $\tau > 0$ and $\phi(t) = \bigl(t/\tau + 1)y_0$, $t \in [-\tau,0]$, for some $y_0 \in \Rset$. For any $t > 0$ the difference between the solution $\hat{y}(t)$ of (\ref{eq:FDDE_Linear_ModifiedOperator})  and the solution $y(t)$ of (\ref{eq:FDDE_Linear}) is
\[
	\hat{y}(t) - y(t) =  \frac{y_0}{\tau} \sum_{k=0}^{\left\lfloor t/\tau \right\rfloor} \lambda^{k} 
		\left[\frac{(t-k\tau)^{\alpha k + 1}}{\Gamma(\alpha k + 2)} - \frac{(t+\tau-k\tau)^{\alpha k + 1}}{\Gamma(\alpha k + 2)} \left( 1 - I_{\frac{\tau}{t+\tau-k\tau}}(2-\alpha, (k + 1)\alpha) \right) \right] .
\]
\end{proposition}

\begin{proof}
Since $\hat{y}(t)$ is solution of (\ref{eq:FDDE_Linear_ModifiedOperator_Equivalent}) as well, we can use Corollary (\ref{cor:FDDE_Lin_Ex_FirstInt_Case2}) to write $\hat{y}(t)$ as solution of (\ref{eq:FDDE_Linear}) with the source term $f(t) + \DC_{0}^{\alpha}  \phi(t) - \DC_{-\tau}^{\alpha} \phi(t)$. Therefore it is 
\begin{equation}\label{eq:ExactSolutionModifiedOperator_1}
	\hat{y}(t) = y(t) + \mathcal{J}_{0,\tau,\lambda} \Bigl( \DC_{0}^{\alpha}  \phi(t) - \DC_{-\tau}^{\alpha} \phi(t) \Bigr) 
\end{equation}
and we can apply Lemma \ref{eq:RL_Integral_PowerFunctions} to evaluate  
\[
	\begin{aligned}
	\mathcal{J}_{0,\tau,\lambda} &\Bigl( \DC_{0}^{\alpha}  \phi(t) - \DC_{-\tau}^{\alpha} \phi(t) \Bigr) = \\
	&=\frac{y_0}{\tau} \sum_{k=0}^{\left\lfloor t/\tau \right\rfloor}  \frac{\lambda^{k}}{\Gamma\bigl(\alpha k + \alpha \bigr)} \int_{0}^{t-k \tau} (t-k \tau - r )^{ \alpha k+ \alpha-1}  \frac{  r^{1-\alpha}-(r+\tau)^{1-\alpha} }{\Gamma(2-\alpha)} \du r = \\
	&= \frac{y_0}{\tau} \sum_{k=0}^{\left\lfloor t/\tau \right\rfloor} \lambda^{k} 
		\left[\frac{(t-k\tau)^{\alpha k + 1}}{\Gamma(\alpha k + 2)} - \frac{(t+\tau-k\tau)^{\alpha k + 1}}{\Gamma(\alpha k + 2)} \left( 1 - I_{\frac{\tau}{t+\tau-k\tau}}(2-\alpha,( k + 1)\alpha) \right) \right] . \\
	\end{aligned}
\] 
which allows to conclude the proof.
\end{proof}

It is also possible to provide the exact solution $\hat{y}(t)$ of (\ref{eq:FDDE_Linear_ModifiedOperator}) in a more compact form.

\begin{proposition}
Let $\tau > 0$ and $\phi(t) = \bigl(t/\tau + 1)y_0$, $t \in [-\tau,0]$, for some $y_0 \in \Rset$. For any $t > 0$ the solution $\hat{y}(t)$ of (\ref{eq:FDDE_Linear_ModifiedOperator})  is
\[
\hat{y}(t)=\frac{y_0}{\tau}\left[\lambda^{p+1}\frac{(t-p\tau)^{\alpha(p+1)+1}}{\Gamma(\alpha(p+1)+2)}+\sum_{k=0}^{p} \lambda^k \frac{(t+\tau-k\tau)^{\alpha k+1}}{\Gamma(\alpha k+2)} I_{\frac{\tau}{t+\tau-k\tau}}(2-\alpha,(k+1)\alpha)\right] + \mathcal{J}_{0,\tau,\lambda} f(t),
\]
where $p\in\mathbb{Z}_+$ is such that $t\in[p\tau,(p+1)\tau)$.
\end{proposition}

\begin{proof}
From Eq. (\ref{eq:FDDE_Lin_Ex_FirstInt_Case2_StepFunction}) and after replacing the source term $f(t)$ with $f(t) + \DC_{0}^{\alpha}  \phi(t) - \DC_{-\tau}^{\alpha} \phi(t)$ we obtain 
\[
	y(t) 
	= \left[ 1 +\frac{1}{\tau} \sum_{k=0}^{\lfloor {t}/{\tau}\rfloor} \lambda^{k+1} \left(u_{k\tau }^{[(k+1)\alpha+1]}(t)-u_{(k+1)\tau}^{[(k+1)\alpha+1]}(t)\right) \right] y_0  
	+\mathcal{J}_{0,\tau,\lambda} f(t) +
	\mathcal{J}_{0,\tau,\lambda} \Bigl( \DC_{0}^{\alpha}  \phi(t) - \DC_{-\tau}^{\alpha} \phi(t) \Bigr) 
\]
and hence, thanks to Proposition  \ref{prop:DifferenceExactSolutions} and by reordering some terms it is  
\begin{align*}
	y(t) 
	&= y_0 +\frac{y_0}{\tau} \sum_{k=0}^{\lfloor{t}/{\tau}\rfloor} \lambda^{k+1} u_{k\tau }^{[(k+1)\alpha+1]}(t)-\frac{y_0}{\tau} \sum_{k=1}^{\lfloor{t}/{\tau}\rfloor+1} \lambda^{k}u_{k\tau}^{[k\alpha+1]}(t) + \mathcal{J}_{0,\tau,\lambda} f(t) \\ 
	&\quad +\frac{y_0}{\tau}\sum_{k=0}^{\lfloor{t}/{\tau}\rfloor} \lambda^k u_{k\tau}^{[k\alpha+1]}(t)-\frac{y_0}{\tau}\sum_{k=-1}^{\lfloor{t}/{\tau}\rfloor-1} \lambda^{k+1} u_{k\tau}^{[(k+1)\alpha+1]}(t)\\
	&\quad+ \frac{y_0}{\tau}\sum_{k=0}^{\lfloor{t}/{\tau}\rfloor} \lambda^k u_{(k-1)\tau}^{[k\alpha+1]}(t) I_{\frac{\tau}{t+\tau-k\tau}}(2-\alpha,(k+1)\alpha)
\end{align*}
We now denote $p=\lfloor{t}/{\tau}\rfloor$ and we simplify some terms in order to obtain 
\[
\begin{aligned}
	y(t) 
	&= y_0 + \frac{y_0}{\tau}\Biggl[\lambda^{p+1}u_{p\tau}^{[(p+1)\alpha+1]}(t) - \lambda^{p+1}u_{(p+1)\tau}^{[(p+1)\alpha]+1}(t)+u_0^{[1]}(t)-u_{-\tau}^{[1]}(t) \Biggr] + \\
	& \quad + \mathcal{J}_{0,\tau,\lambda} f(t) 
	+\frac{y_0}{\tau} \sum_{k=0}^{p} \lambda^k u_{(k-1)\tau}^{[k\alpha+1]}(t) I_{\frac{\tau}{t+\tau-k\tau}}(2-\alpha,(k+1)\alpha).
\end{aligned}
\]
Observing that $u_{(p+1)\tau}^{[(p+1)\alpha]+1}(t)=0$, as $t\in[p\tau,(p+1)\tau)$, and $u_0^{[1]}(t)-u_{-\tau}^{[1]}(t)=-\tau$, we finally obtain
\[
\begin{aligned}
	y(t) 
	&= \frac{y_0}{\tau}\left[ \lambda^{p+1}u_{p\tau}^{[(p+1)\alpha+1]}(t)
	+ \sum_{k=0}^{p} \lambda^k u_{(k-1)\tau}^{[k\alpha+1]}(t) I_{\frac{\tau}{t+\tau-k\tau}}(2-\alpha,(k+1)\alpha)\right] + \mathcal{J}_{0,\tau,\lambda} f(t) 
\end{aligned}
\]
from which the proof follows.
\end{proof}

\begin{figure}[tph]
\centering
\includegraphics[width=0.95\textwidth]{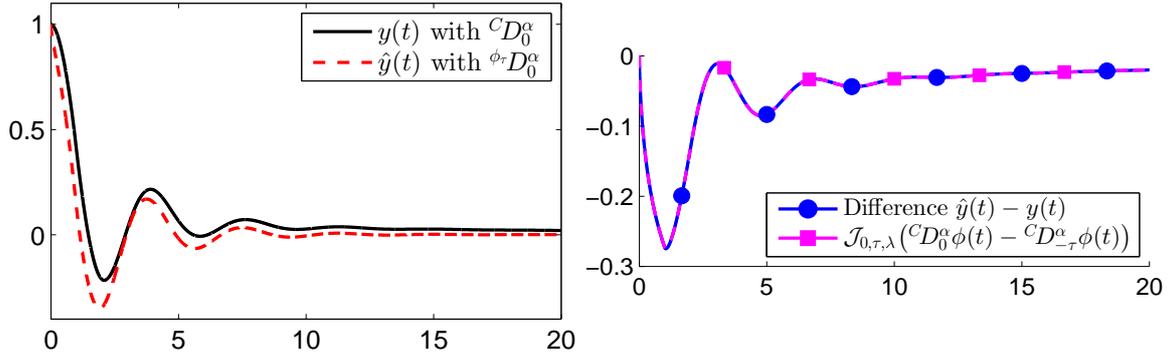}
\caption{Solutions of the linear FDDEs (\ref{eq:FDDE_Linear}) and (\ref{eq:FDDE_Linear_ModifiedOperator}) for $\lambda=-1.0$, $\tau=1.0$, $f(t)=0$ and $\phi(t)=\bigl(t/\tau + 1)y_0$}\label{Fig_FDDE_Test_Linear_2}
\end{figure}

In order to show the different behaviors due to the two different operators $\DC^{\alpha}_0$ and $\Dphi^{\alpha}_0$, in the left plot of Figure \ref{Fig_FDDE_Test_Linear_2} we present the solutions $y(t)$ and $\hat{y}(t)$ of the linear FDDEs (\ref{eq:FDDE_Linear}) and (\ref{eq:FDDE_Linear_ModifiedOperator}). We consider here a problem without forcing term (namely $f(t)\equiv 0$). The difference $\hat{y}(t)-y(t)$ between the two solutions is presented in the right plot together with the generalized integral $\mathcal{J}_{0,\tau,\lambda} \bigl( \DC_{0}^{\alpha}  \phi(t) - \DC_{-\tau}^{\alpha} \phi(t) \bigr)$; one can clearly appreciate that the two plots overlap, thus confirming that $\hat{y}(t)-y(t) = \mathcal{J}_{0,\tau,\lambda} \bigl( \DC_{0}^{\alpha}  \phi(t) - \DC_{-\tau}^{\alpha} \phi(t) \bigr)$ as expected from Eq. (\ref{eq:ExactSolutionModifiedOperator_1}).

We observe that the main difference between the two solutions occurs for small values of time $t$, whilst asymptotically they tend to coincide as $t\rightarrow\infty$.

\begin{figure}[tph]
\centering
\includegraphics[width=0.95\textwidth]{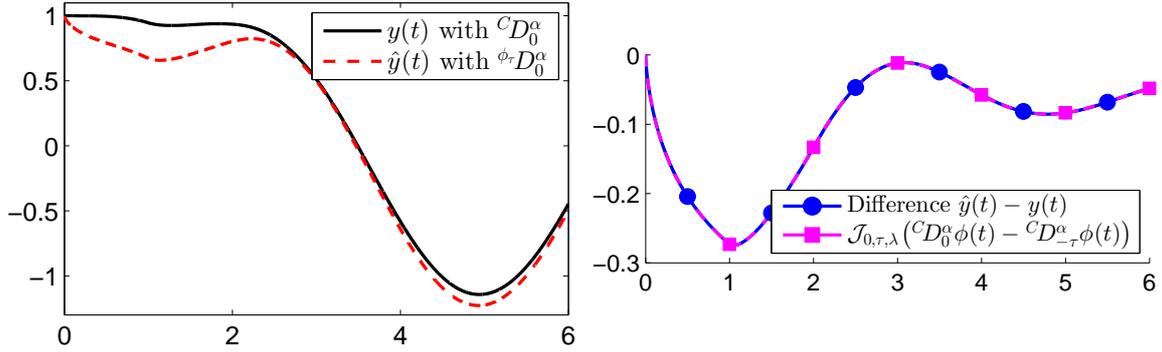}
\caption{Solutions of the linear FDDEs (\ref{eq:FDDE_Linear}) and (\ref{eq:FDDE_Linear_ModifiedOperator}) for $\lambda=-1.0$, $\tau=1.0$, $f(t)=\sin t$ and $\phi(t)=\bigl(t/\tau + 1)y_0$}\label{Fig_FDDE_Test_Linear_3}
\end{figure}

Similar results are presented in Figure \ref{Fig_FDDE_Test_Linear_3} for a problem with a source term $f(t)=\sin t$. 


\section{Nonlinear FDDEs: solution by means of numerical methods}\label{S:Nonlinear}

Finding explicit solutions of the nonlinear FDDE (\ref{eq:FDDE_General}) is, in general, not possible and therefore, numerical methods are necessary.

Since the aim of this paper is just to highlight the impact of the initial data on the solution and on the fractional operator, and not to devise highly efficient numerical methods, we consider here basic methods with just a first-order accuracy with respect to the step-size. We refer to the existing literature for specific works concerning numerical methods for FDDEs (e.g., see 
\cite{morgado2013analysis,CaoLuo2018,DabiriButcher2018,HendyPimenovMaciasDiaz2020,JhingaDaftardarGejji2019}).

As usual, on the integration interval $[0,T]$ we consider an equispaced grid $t_n=nh$, $n=0,1,\dots,N$, where $N=\left\lceil T/h\right\rceil$ and $h>0$ is the step-size. 

To compute numerical approximations of the solution of FDDEs with the standard Caputo derivative $\DC^{\alpha}_{0}$ we can apply a standard product-integration rule based on the integral formulation of (\ref{eq:FDDE_General}) which, for $0<\alpha<1$, reads as 
\[
	y(t_n) = y_0 + \frac{1}{\Gamma(\alpha)} \sum_{j=0}^{n-1}\int_{t_j}^{t_{j+1}} (t_n-u)^{\alpha-1} g(u,y(u),y(u-\tau)) \du u ,
\]
and where in each interval $[t_j,t_{j+1}]$ the vector field $g(u,y(u),y(u-\tau))$ is approximated by the constant $g(t_j,y_j,y(t_j - \tau))$. The numerical approximation is hence given by
\[
	y_n = y_0 + h^{\alpha} \sum_{j=0}^{n-1} b_{n-j} g(t_j,y_j,y(t_j - \tau))
	, \quad
	b_n = \frac{n^{\alpha} - (n-1)^{\alpha}}{\Gamma(\alpha+1)} .
\]

This method is known as the rectangular product-integration rule or $1$-step Adams-Bashforth method and it is widely used and studied in fractional calculus since the pioneering works by Diethelm and co-authors \cite{DiethelmFordFreed2002,DiethelmFordFreed2004}.

The simplest approximation scheme for the operator $\Dphi_{{\scriptscriptstyle 0}}^{\alpha}$ is instead obtained by just fixing $h>0$ in (\ref{eq:FracDerviPhi_CalcualtingFormula}) and replacing the finite differences in (\ref{eq:FDDE_General}) thus to obtain the computational scheme 
\[
	y_n = \phi(-\tau) - \sum_{j=1}^{J_1} \omega_{j}^{(\alpha)} \bigl( y_{n-j} - \phi(-\tau) \bigr) - \sum_{j=J_1+1}^{J_2} \omega_{j}^{(\alpha)} \bigl( \phi(t_n-jh) - \phi(-\tau) \bigr) + h^{\alpha} g(t_n,y_n,y(t_n - \tau)) ,
\]
where $J_1 = \left\lfloor t_n/h \right\rfloor$ and $J_2 = \left\lfloor (t_n+\tau)/h \right\rfloor$. An error ${\mathcal O}\bigl(h\bigr)$, $h\to 0$, is expected as in the usual GL scheme derived in the same way from (\ref{eq:GL_Definition}). The last term $g(t_n,y_n,y(t_n - \tau))$ can be hence replaced by $g(t_{n-1},y_{n-1},y(t_{n-1} - \tau))$ in order to obtain an explicit scheme and reduce the computational complexity.

Note that in both schemes the values $y(t_j - \tau)$ are not available if $t_j-\tau$ is outside the mesh-grid $\bigl\{ t_n\bigr\}_{n\in\Nset}$. Therefore it can be necessary to select the two closest mesh points $t_{k}$ and $t_{k+1}$, i.e. such that $t_j-\tau \in [t_{k} , t_{k+1}]$, and perform an interpolation of $(t_k,y_k)$ and $(t_{k+1},y_{k+1})$ to obtain a suitable approximation of $y(t_j - \tau)$. A first-order approximation is however sufficient since both methods are just first-order convergent. The use of interpolation avoids constraining the step-size $h$ to the delay $\tau$ since the method can operate also when $t_j-\tau$ is outside the mesh. In the numerical simulation we use a step-size $h=2^{-8} \approx 3.9 \times 10^{-3}$; since the error is ${\mathcal O}\bigl(h\bigr)$ this is sufficient to provide an accurate enough approximation for graphical visualization.   

\begin{figure}[tph]
\centering
\includegraphics[width=0.95\textwidth]{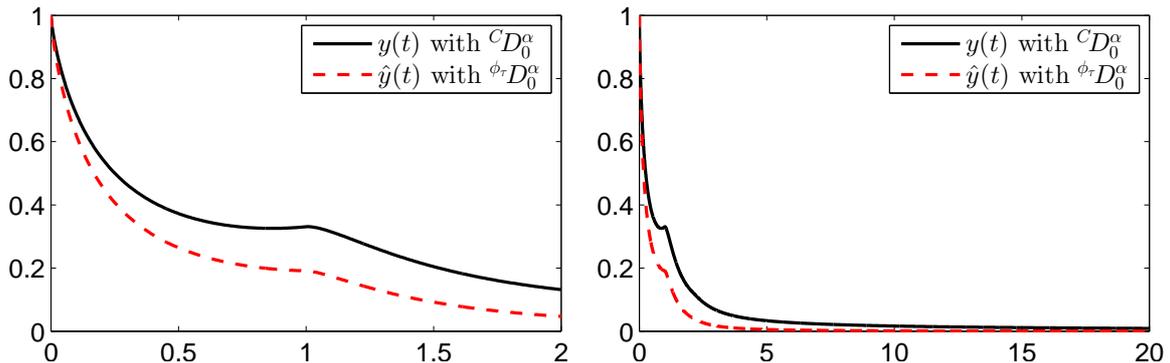}
\caption{Short and long-term solutions of the linear FDDEs (\ref{eq:FDDE_Linear}) and (\ref{eq:FDDE_Linear_ModifiedOperator}) for $\alpha=0.8$, $g(t,y(t),y(t-\tau)) =  -2y(t)(1.2-y(t-\tau)$, $\tau=1.0$, and $\phi(t)=\bigl(t/\tau + 1)y_0$.}\label{fig:Fig_FDDE_Test_NonLinear_1}
\end{figure}

In Figure \ref{fig:Fig_FDDE_Test_NonLinear_1} we present the results of the simulations for a nonlinear problem where $g(t,y(t),y(t-\tau)) =  -2y(t)\left(1.2-y(t-\tau)\right)$ and the  initial condition $\phi(t)=\left(t/\tau + 1\right)y_0$, $-\tau\le t \le 0$; the left plot shows the differences between the solutions obtained with the operators $\DC^{\alpha}_0$ and $\Dphi^{\alpha}_0$, respectively close to the origin, while the right plot highlights how the two solutions tend to overlap over long-time integration.

\begin{figure}[tph]
\centering
\includegraphics[width=0.95\textwidth]{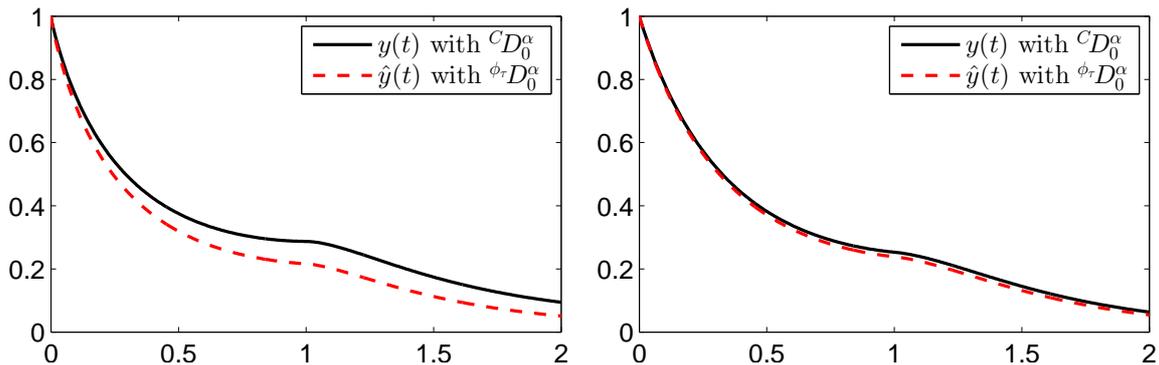}
\caption{Solutions of the linear FDDEs (\ref{eq:FDDE_Linear}) and (\ref{eq:FDDE_Linear_ModifiedOperator}) for $g(t,y(t),y(t-\tau)) =  -2y(t)(1.2-y(t-\tau)$, $\tau=1.0$, and $\phi(t)=\bigl(t/\tau + 1)y_0$ and $\alpha=0.9$ (left plot) and $\alpha=0.98$ (right plot).}\label{fig:Fig_FDDE_Test_NonLinear_2}
\end{figure}

For the same test problem we illustrate the results when the order $\alpha$ approaches 1. In the left plot of Figure \ref{fig:Fig_FDDE_Test_NonLinear_2} we used $\alpha=0.9$ while in the right plot we used $\alpha=0.98$. As expected, the two solutions tend to coincide as $\alpha \to 1$; the difference between the operators $\DC^{\alpha}_0$ and $\Dphi^{\alpha}_0$ tends to vanish in this case, as the phenomenon under investigation is a peculiarity of FDDEs and not of integer-order DDEs.

\section{Concluding remarks}\label{S:ConcludingRemarks}

We have discussed some issues related to the initialization of FDDEs. We have observed that when the fractional derivative is intended as a generalization (by means of the GL process) of the usual limit of the difference quotient defining integer-order derivatives, the function must be forced to assume specific values at the left of the initial point in order to retrieve the Caputo derivative. Hence, initial conditions of FDDEs assign values of the solution before the initial point and when these values do not match with the values forced to obtain the Caputo derivative, some inconsistencies arise.

For this reason, we think that it is not advisable to use the usual Caputo derivative in FDDEs and we have proposed a fractional operator obtained on the basis of the GL derivative, suitably modified in order to match the initial condition of the FDDE. When the initial condition $\phi(t)$ is equal to the constant value $y_0$ assumed at the left of the origin to obtain the Caputo derivative, the two operators coincide. Otherwise, this process leads to a different operator which is studied in this paper. It is beyond the scope of this paper to introduce a further fractional derivative; we have just included the initial condition of the FDDE in the process leading to the fractional derivative, in order to avoid inconsistencies.

The above reasoning does not imply that it is wrong or unacceptable to incorporate the usual Caputo derivative in FDDEs,  based on the left-inverse of the RL integral. Nevertheless, in this case, the Caputo derivative cannot any longer be considered as a generalization of the integer-order derivative, but just as one of the left-inverses of the RL integral and, perhaps, the term derivative should be avoided. Alternatively, one could restrict to use only initial functions matching the assumptions which are necessary to obtain the equivalence between the derivative based on the GL definition and the Caputo derivative; however, such kind of a limitation could turn out to be overly  restrictive in several models.

Moreover, this paper does not have the ambition of providing definitive answers, but only to highlight possible inconsistencies in the definition of problems with fractional derivatives and encourage more in-depth investigations in this field.

For shortness, we have not discussed FDDEs with the RL or other fractional-order derivatives, but the analysis can be surely extended to cover other cases as well.

\appendix
\section{Alternative proof}

The exact solution of linear FDDEs (\ref{eq:FDDE_Linear}) has been found in Proposition \ref{prop:DDE_Frac_ExSol}, as well as Corollaries \ref{cor:FDDE_Lin_Ex_FirstInt_Case1} and \ref{cor:FDDE_Lin_Ex_FirstInt_Case2}, by using the LT under suitable assumptions for the exponential boundedness of the exact solution. Although the LT is useful for deriving the exact solution of linear FDDEs depending on the initial function $\phi(t)$, the results presented in Section \ref{S:ExactSolution} are more general and can be proved without employing the LT and the assumptions necessary for using this tool, such as inequality \eqref{eq:AssumptionLinearExpBound}.

We present here, just for completeness, a proof of Corollary \ref{cor:FDDE_Lin_Ex_FirstInt_Case1} based on standard mathematical induction arguments. For shortness, we only consider the solution obtained for the initial condition $\phi(t)=y_0$, $-\tau \le t \le 0$. Further initial conditions can be treated in similar way.

Just for notational convenience we reformulate the statement of Corollary \ref{cor:FDDE_Lin_Ex_FirstInt_Case1} in a slightly different form.

\begin{corollary}
Let $\tau > 0$ and $\phi(t)\equiv y_0$, $t \in [-\tau,0]$, for some $y_0 \in \Rset$. For any $t \ge 0$ the exact solution of the linear FDDE \eqref{eq:FDDE_Linear} is
\begin{align}\label{exact.sol}
    y(t)&=y_0\sum_{k=0}^\infty\lambda^{k}u_{(k-1)\tau}^{[k\alpha]}(t)+\sum_{k=0}^\infty\lambda^k\int_{0}^t u_{k\tau}^{[(k+1)\alpha-1]}(t-r)f(r)dr\\
 \nonumber   &=y_0\sum_{k=0}^{\lfloor\frac{t}{\tau}\rfloor+1}\lambda^{k}u_{(k-1)\tau}^{[k\alpha]}(t)+\sum_{k=0}^{\lfloor\frac{t}{\tau}\rfloor}\lambda^{k}(u_{k\tau}^{[(k+1)\alpha-1]}\ast f)(t)
\end{align}
\end{corollary}

\begin{proof}
We proceed by mathematical induction. The first step is to show that \eqref{exact.sol} holds whenever $t\in[0,\tau)$. This follows in a straightforward way, by a direct application of the RL integral operator $J_0^\alpha$ to both sides of the fractional-order differential equation  \eqref{eq:FDDE_Linear}, considering $t\in[0,\tau)$. 

For the second induction step,
let $p\in\mathbb{Z}_+$, $p\geq 1$, and let us assume that the formula \eqref{exact.sol} holds for any $t\in[0, p\tau)$. It remains to show that \eqref{exact.sol} is true for any $t\in[p\tau,(p+1)\tau)$. Indeed, starting from the FDDE (\ref{eq:FDDE_Linear}) and applying the RL integral operator $J_0^\alpha$ and considering $t\in[p\tau,(p+1)\tau)$, it follows that: 
$$y(t)-y_0=\frac{\lambda}{\Gamma(\alpha)}\int_0^t (t-r)^{\alpha-1}y(r-\tau)dr+\frac{1}{\Gamma(\alpha)}\int_0^t (t-r)^{\alpha-1} f(r)dr$$
and therefore, by the induction hypothesis: 
\begin{align*}
    y(t)&=y_0+\frac{\lambda}{\Gamma(\alpha)}\left[\sum_{j=0}^{p-1}\int_{j\tau}^{(j+1)\tau}(t-r)^{\alpha-1}y(r-\tau)dr+\int_{p\tau}^t (t-r)^{\alpha-1}y(r-\tau)dr\right]+(u_0^{[\alpha-1]}\ast f)(t)\\
    &=y_0+\frac{\lambda}{\Gamma(\alpha)}\left[\int_0^\tau(t-r)^{\alpha-1}dr\right]y_0+\\
    &\quad +\frac{\lambda}{\Gamma(\alpha)} \sum_{j=1}^{p-1}\int_{j\tau}^{(j+1)\tau}(t-r)^{\alpha-1}\left[\sum_{k=0}^{j}\lambda^{k}u_{(k-1)\tau}^{[k\alpha]}(r-\tau)y_0+\sum_{k=0}^{j-1}\lambda^{k}(u_{k\tau}^{[(k+1)\alpha-1]}\ast f)(r-\tau)\right]dr +\\
    &\quad +\frac{\lambda}{\Gamma(\alpha)}\int_{p\tau}^t (t-r)^{\alpha-1}\left[\sum_{k=0}^{p}\lambda^{k}u_{(k-1)\tau}^{[k\alpha]}(r-\tau)y_0+\sum_{k=0}^{p-1}\lambda^k(u_{k\tau}^{[(k+1)\alpha-1]}\ast f)(r-\tau)\right]dr+\\
    &\quad +(u_0^{[\alpha-1]}\ast f)(t).
    \end{align*}
Grouping the terms conveniently, leads to:
\begin{align*}
  y(t)  &=y_0+\frac{\lambda}{\Gamma(\alpha)}\left[\sum_{j=0}^{p-1}\int_{j\tau}^{(j+1)\tau}(t-r)^{\alpha-1}dr+\int_{p\tau}^t (t-r)^{\alpha-1}dr\right]y_0+\\
    &\quad + \frac{\lambda}{\Gamma(\alpha)}\left[\sum_{j=1}^{p-1}\sum_{k=1}^{j}\lambda^{k}\int_{j\tau}^{(j+1)\tau}(t-r)^{\alpha-1}u_{(k-1)\tau}^{[k\alpha]}(r-\tau)dr+\sum_{k=1}^{p}\lambda^{k}\int_{p\tau}^t(t-r)^{\alpha-1}u_{(k-1)\tau}^{[k\alpha]}(r-\tau)dr\right] y_0+\\
     &\quad + \frac{\lambda}{\Gamma(\alpha)}\sum_{j=1}^{p-1}\sum_{k=0}^{j-1}\lambda^{k}\int_{j\tau}^{(j+1)\tau}(t-r)^{\alpha-1}(u_{k\tau}^{[(k+1)\alpha-1]}\ast f)(r-\tau)dr+\\
     &\quad + \frac{\lambda}{\Gamma(\alpha)}\sum_{k=0}^{p-1}\lambda^{k}\int_{p\tau}^t(t-r)^{\alpha-1}(u_{k\tau}^{[(k+1)\alpha-1]}\ast f)(r-\tau)dr +(u_0^{[\alpha-1]}\ast f)(t).
\end{align*}
Index changes in the double summations lead to:
\begin{align*}
  y(t)  &=\left[1+\frac{\lambda}{\Gamma(\alpha)}\int_{0}^{t}(t-r)^{\alpha-1}dr\right]y_0+\\
   &\quad + \frac{\lambda}{\Gamma(\alpha)}\left[\sum_{k=1}^{p-1}\sum_{j=k}^{p-1}\lambda^{k}\int_{j\tau}^{(j+1)\tau}(t-r)^{\alpha-1}u_{(k-1)\tau}^{[k\alpha]}(r-\tau)dr+\sum_{k=1}^{p}\lambda^{k}\int_{p\tau}^t(t-r)^{\alpha-1}u_{(k-1)\tau}^{[k\alpha]}(r-\tau)dr\right] y_0+\\
     &\quad + \frac{\lambda}{\Gamma(\alpha)}\sum_{k=0}^{p-2}\sum_{j=k+1}^{p-1}\lambda^{k}\int_{j\tau}^{(j+1)\tau}(t-r)^{\alpha-1}(u_{k\tau}^{[(k+1)\alpha-1]}\ast f)(r-\tau)dr+\\
     &\quad + \frac{\lambda}{\Gamma(\alpha)}\sum_{k=0}^{p-1}\lambda^{k}\int_{p\tau}^t(t-r)^{\alpha-1}(u_{k\tau}^{[(k+1)\alpha-1]}\ast f)(r-\tau)dr+  +(u_0^{[\alpha-1]}\ast f)(t).
\end{align*}
It is easy to see that the above formula can be simplified to:
\begin{align*}
  y(t)  &=\left[1+\lambda u_0^{[\alpha]}(t) + \frac{1}{\Gamma(\alpha)}\sum_{k=1}^{p}\lambda^{k+1}\int_{k\tau}^{t}(t-r)^{\alpha-1}u_{(k-1)\tau}^{[k\alpha]}(r-\tau)dr\right] y_0+\\
     &\quad +(u_0^{[\alpha-1]}\ast f)(t)+ \frac{1}{\Gamma(\alpha)}\sum_{k=1}^{p}\lambda^{k}\int_{k\tau}^{t}(t-r)^{\alpha-1}(u_{(k-1)\tau}^{[k\alpha-1]}\ast f)(r-\tau)dr.
\end{align*}
Using Lemma \ref{eq:RL_Integral_PowerFunctions}, we evaluate the integrals above: 
\begin{align*}\frac{1}{\Gamma(\alpha)}&\int_{k\tau}^t (t-r)^{\alpha-1}  u_{(k-1)\tau}^{[k\alpha]}(r-\tau)dr=\frac{1}{\Gamma(\alpha)\Gamma(k\alpha+1)}\int_{k\tau}^t (t-r)^{\alpha-1}(r-k\tau)^{k\alpha}dr=\\
&=\frac{1}{\Gamma(\alpha)\Gamma(k\alpha+1)}\int_{0}^{t-k\tau} (t-k\tau-r)^{\alpha-1}r^{k\alpha}dr=\frac{(t-k\tau)^{(k+1)\alpha}}{\Gamma((k+1)\alpha+1)}=u_{k\tau}^{[(k+1)\alpha]}(t).
\end{align*}
On the other hand, with a change in the order of integration and making use of Lemma \ref{eq:RL_Integral_PowerFunctions} again, we obtain: 
\begin{align*}
\frac{1}{\Gamma(\alpha)}&\int_{k\tau}^{t}(t-r)^{\alpha-1} (u_{(k-1)\tau}^{[k\alpha-1]}\ast f)(r-\tau)dr= \\
&=\frac{1}{\Gamma(\alpha)\Gamma(k\alpha)}\int_{k\tau}^t(t-r)^{\alpha-1}\left(\int_{0}^{r-k\tau}(r-k\tau-s)^{k\alpha-1}f(s)ds\right)dr\\
&=\frac{1}{\Gamma(\alpha)\Gamma(k\alpha)}\int_{0}^{t-k\tau}f(s)\left(\int_{s+k\tau}^{t}(t-r)^{\alpha-1}(r-k\tau-s)^{k\alpha-1}dr\right)ds\\
&=\frac{1}{\Gamma(\alpha)\Gamma(k\alpha)}\int_{0}^{t-k\tau}f(s)\left(\int_{0}^{t-k\tau-s}(t-k\tau-s-r)^{\alpha-1}r^{k\alpha-1}dr\right)ds\\
&=\int_{0}^{t-k\tau}f(s)\frac{(t-k\tau-s)^{(k+1)\alpha-1}}{\Gamma((k+1)\alpha)}ds\\
&=\int_{0}^{t-k\tau}f(s)u_{k\tau}^{[(k+1)\alpha-1]}(t-s)ds\\
&=\int_{0}^{t}f(s)u_{k\tau}^{[(k+1)\alpha-1]}(t-s)ds
\\&=
(u_{k\tau}^{[(k+1)\alpha-1]}\ast f)(t).
\end{align*}
We finally obtain:
$$
  y(t)  =y_0\sum_{k=0}^{p+1}\lambda^{k}u_{(k-1)\tau}^{[k\alpha]}(t)+ \sum_{k=0}^{p}\lambda^{k}(u_{k\tau}^{[(k+1)\alpha-1]}\ast f)(t),\quad\forall~t\in[p\tau,(p+1)\tau),
$$
which completes the proof.
\end{proof}

\bibliographystyle{elsarticle-num}
\bibliography{FDDE_IC_Biblio}

\end{document}